\documentclass[12pt,reqno]{amsart}
\usepackage{amsfonts}
\usepackage{amssymb, amsmath, latexsym}
\usepackage{eufrak}
\usepackage{lipsum} % for filler text
\setbox0=\hbox{$+$}
\newdimen\plusheight
\plusheight=\ht0
\def\+{\;\lower\plusheight\hbox{$+$}\;}

\setbox0=\hbox{$-$}
\newdimen\minusheight
\minusheight=\ht0
\def\-{\;\lower\minusheight\hbox{$-$}\;}

\setbox0=\hbox{$\cdots$}
\newdimen\cdotsheight
\cdotsheight=\plusheight%\ht0
\def\cds{\lower\cdotsheight\hbox{$\cdots$}}

\renewcommand{\(}{\left\(}
\renewcommand{\)}{\right\)}
\renewcommand{\[}{\left[}

\numberwithin{equation}{section}
 \theoremstyle{plain}
\newtheorem{theorem}{Theorem}[section]
\newtheorem{lemma}[theorem]{Lemma}
\newtheorem{conjecture}[theorem]{Conjecture}
\newtheorem{corollary}[theorem]{Corollary}

\newenvironment{pf-thm-nu-50n8}
   {\vskip 0.15in \par\noindent{\it Proof of Theorem \ref{thm-nu-50n8}.}\hskip 0.5em\ignorespaces}
    {\hfill $\Box$\par\medskip}

\newenvironment{pf-cor1-nu-50n8}
   {\vskip 0.15in \par\noindent{\it Proof of \eqref{p-nu-250n208}.}\hskip 0.5em\ignorespaces}
    {\hfill $\Box$\par\medskip}

\newenvironment{pf-cor2-nu-50n8}
   {\vskip 0.15in \par\noindent{\it Proof of Corollary \ref{cor-nu-50n8}.}\hskip 0.5em\ignorespaces}
    {\hfill $\Box$\par\medskip}

 \newenvironment{pf-thm-spt-main}
   {\vskip 0.15in \par\noindent{\it Proof of Theorem \ref{thm-spt-main}.}\hskip 0.5em\ignorespaces}
    {\hfill $\Box$\par\medskip}

 \newenvironment{pf-spt-1250n573}
   {\vskip 0.15in \par\noindent{\it Proofs of \eqref{cong250n73} and \eqref{cong1250n573}.}\hskip 0.5em\ignorespaces}
    {\hfill $\Box$\par\medskip}

  \newenvironment{pf-thm-sptbar}
   {\vskip 0.15in \par\noindent{\it Proof of Theorem \ref{thm-sptbar}.}\hskip 0.5em\ignorespaces}
    {\hfill $\Box$\par\medskip}

     \newenvironment{pf-cor-sptbar}
   {\vskip 0.15in \par\noindent{\it Proof of Corollary \ref{cor-sptbar}.}\hskip 0.5em\ignorespaces}
    {\hfill $\Box$\par\medskip}

\newenvironment{pf-thm-p-omega}
   {\vskip 0.15in \par\noindent{\it Proof of Theorem  \ref{thm-p-omega}.}\hskip 0.5em\ignorespaces}
    {\hfill $\Box$\par\medskip}

    \newenvironment{pf-cor-p-omega}
   {\vskip 0.15in \par\noindent{\it Proof of Corollary \ref{cor-p-omega}.}\hskip 0.5em\ignorespaces}
    {\hfill $\Box$\par\medskip}

\begin{document}
\allowdisplaybreaks
\title[Generating Functions and Congruences for Some Partition Functions Related to Mock Theta Functions] {Generating Functions and Congruences for Some Partition Functions Related to Mock Theta Functions}

\author{Nayandeep Deka Baruah}
\address{Department of Mathematical Sciences, Tezpur University, Sonitpur, Assam, India, Pin-784028}
\email{nayan@tezu.ernet.in}

\author{Nilufar Mana Begum}
\address{Department of Mathematical Sciences, Tezpur University, Sonitpur, Assam, India, Pin-784028}
\email{nilufar@tezu.ernet.in}

%\maketitle

\begin{center}
{\textbf{Generating Functions and Congruences for Some Partition Functions Related to Mock Theta Functions}}\\[5mm]
{\footnotesize  Nayandeep Deka Baruah and Nilufar Mana Begum}\\[3mm]
\end{center}

\vskip 5mm \noindent{\footnotesize{\bf Abstract.}
 Recently, Andrews, Dixit and Yee introduced partition functions associated with Ramanujan/Watson third order mock theta functions $\omega(q)$ and $\nu(q)$. In this paper, we find several new exact generating functions for those partition functions as well as the associated smallest parts functions and deduce several new congruences modulo powers of 5.

\vskip 3mm
\noindent {\footnotesize 2010 Mathematical Reviews Classification
Numbers: Primary 11P83; Secondary 05A15, 05A17}.}

\section{\textbf{Introduction}}

A \emph{partition}
$\lambda=(\lambda_1,\lambda_2,\cdots,\lambda_k)$ of a positive integer
$n$ is a finite sequence of non-increasing positive integer
\emph{parts} $\lambda_i$ such that $n=\sum_{i=1}^k\lambda_i$. For example, the partitions of $5$ are $(5), ~(4,1),~ (3,2),~
(3,1,1),~ (2,2,1),~ (2,1,1,1),~ (1,1,1,1,1)$.

Let $p(n)$ denote the number of partitions of a positive integer $n$. With the convention that $p(0)=1$, the generating
function for $p(n)$ is given by
\begin{align}\label{p-n}
\sum_{n=0}^\infty p(n)q^n=\frac{1}{(q;q)_\infty},
\end{align}
where, here and for the sequel, we use the customary $q$-series notation:
\begin{align*}(a;q)_0&:=1,\quad (a;q)_n:=\prod_{k=0}^{n-1}(1-aq^k), \quad n\ge 1,\\\intertext{and}
(a;q)_\infty&:=\lim_{n\to \infty}(a;q)_n, \quad |q|<1.\end{align*} For any positive integer $j$, for brevity, we also use $E_j:=(q^j;q^j)_\infty$.

In 1919, Ramanujan \cite{rama-pn} proved that
\begin{align}\label{rama-beautiful}
\sum_{n=0}^\infty p(5n+4)q^n=5~\dfrac{E_5^5}{E_1^6},
\end{align}
which immediately implies one of his three famous partition congruences, namely,
\begin{align*}
p(5n+4)&\equiv0~(\textup{mod}~5).
\end{align*}

In \cite{baruah-begum}, the authors found  exact generating functions for $Q(5n+1)$, $Q(25n+1)$ and $Q(125n+26)$, where $Q(n)$ denotes the number of partitions of a nonnegative integer into  distinct (or, odd) parts. In sequel, in this paper, we find several new exact generating functions for some partition functions associated with Ramanujan/Watson third order mock theta functions as well as the associated smallest parts functions and deduce several new congruences modulo powers of 5.

Recently, partition-theoretic interpretations of mock theta functions have been the subject of prominent study. Garthwaite \cite{garth} showed the existence of infinitely many congruences for the third order mock theta function
\begin{align}\label{wq}\omega(q):=\sum_{n=0}^\infty \dfrac{q^{2n^2+2n}}{(q;q^2)_{n+1}^2}=\sum_{n=0}^\infty a_\omega(n)q^n.\end{align} However, the first explicit congruences were given by Waldherr \cite{waldherr}:
\begin{align}\label{waldhrr-cong}
a_\omega(40n+27)\equiv a_\omega(40n+35)\equiv 0~(\textup{mod}~5).
\end{align}

Recently, Andrews, Dixit and Yee \cite{andrews-research} introduced partition functions associated with $\omega(q)$ and $\nu(q)$, where the latter one is a
third-order mock theta function,
 $$\nu(q):=\sum_{n=0}^\infty \dfrac{q^{n(n+1)}}{(-q;q^2)_{n+1}}.$$It is worthwhile to note that $\omega(q)$ and $\nu(q)$ are related by \cite[p. 62, Eq. (26.88)]{fine}:
 $$\nu(-q)=q\omega(q^2)+\dfrac{E_4^3}{E_2^2}.$$Let $p_\omega(n)$ denote the number of partitions of $n$ in which each odd part is less than twice the smallest part and let $p_\nu(n)$ denote the number of partitions of $n$ in which the parts are distinct  and all odd parts are less than twice the smallest part. It was shown  by Andrews, Dixit and Yee \cite{andrews-research} that
\begin{align}\label{qwq}\sum_{n=1}^\infty p_\omega(n)q^n=q\omega(q)\\\intertext{and}
\sum_{n=1}^\infty p_\nu(n)q^n=\nu(-q).\end{align} By \eqref{wq} and \eqref{qwq}, it is clear that $p_\omega(n)=a_\omega(n-1)$, and hence, Wladherr's congruences \eqref{waldhrr-cong} can be recast as
 \begin{align}\label{waldherr-cong}
p_\omega(40n+28)\equiv p_\omega(40n+36)\equiv 0~(\textup{mod}~5).
\end{align}Recently, Andrews, Passary, Sellers and Yee \cite{andrews-rama} found an elementary proof of the above congruences. They also proved several congruences modulo 2 and infinite families of congruences modulo 4 and modulo 8 for $p_\omega(n)$ and $p_\nu(n)$. Motivated by the works in \cite{andrews-research, andrews-rama}, Wang \cite{wang} and Cui, Gu and Hao \cite{cui-gu-hao} found many new congruences satisfied by $p_\omega(n)$ and $p_\nu(n)$ modulo 11 and modulo powers of 2 and 3. In particular, Wang \cite{wang} derived the following exact generating functions:
\begin{align}
\label{gen-nu-2n}\sum_{n=0}^\infty p_\nu(2n)q^n&=\dfrac{E_2^{3}}{E_1^2},\\
\label{gen-om-8n4}\sum_{n=0}^\infty p_\omega(8n+4)q^n&=4\dfrac{E_2^{10}}{E_1^9}.
\end{align}
In fact, the first identity was proved earlier by Hirschhorn and Sellers \cite[Eq. (9)]{mike-sellers-shell} while proving some results for the
so-called 1-shell totally symmetric plane partitions, first introduced by Blecher \cite{blecher}. It is to be noted that
\begin{align}\label{nu-1-shell}p_\nu(2n)=f(6n+1),\end{align}  where $f(n)$ counts the number of 1-shell totally symmetric plane partitions of $n$. We refer to \cite{blecher} and \cite{mike-sellers-shell} for definition and further details. Xia \cite{xia-shell} proved that
\begin{align}\label{gen-nu-10n8}\sum_{n=0}^\infty f(30n+25)q^n=\sum_{n=0}^\infty p_\nu(10n+8)q^n=5\dfrac{E_2^{2}E_5^{2}E_{10}}{E_1^4},\end{align}
which is the exact generating function of $p_\nu(10n+8)$ and the congruence \eqref{spt-nu-10n8} below, proved by Andrews, Dixit and Yee \cite{andrews-research}, follows immediately. Xia \cite{xia-shell} also proved that
\begin{align*}
f(750n+625)  \equiv  0~(\textup{mod}~25),
\end{align*}
which is clearly equivalent to
\begin{align}\label{p-nu-250n208}
p_\nu(250n+208)\equiv  0~(\textup{mod}~25).
\end{align}
In this paper, we find the following representation of the generating function of \linebreak
$p_\nu(50n+8)$.

\begin{theorem}\label{thm-nu-50n8}
We have
\begin{align}\label{gen-nu-50n8-main}\sum_{n=0}^\infty p_\nu(50n+8)q^n&=5\Bigg(\dfrac{E_{2}^{5}E_5^4}{E_1^{6}E_{10}^2}+160q\dfrac{E_2^{11}E_5^4}{E_1^{14}}+2000q^2
\dfrac{E_2^{11}E_5^{10}}{E_1^{20}}\Bigg).
\end{align}
\end{theorem}
We also deduce \eqref{p-nu-250n208} and the following new congruence modulo 125.

\begin{corollary}\label{cor-nu-50n8} For any nonnegative integer $n$, we have
\begin{align}\label{p-nu-6250n5208}p_\nu(6250n+5208)&\equiv0~(\textup{mod}~125).
\end{align}
\end{corollary}

From \eqref{nu-1-shell}, Theorem \ref{thm-nu-50n8} and Corollary \ref{cor-nu-50n8}, we also have the following new results on $f(n)$, the number of 1-shell totally symmetric plane partitions of $n$.

\begin{theorem}
We have
\begin{align*}\sum_{n=0}^\infty f(150n+25)q^n&=5\Bigg(\dfrac{E_{2}^{5}E_5^4}{E_1^{6}E_{10}^2}+160q\dfrac{E_2^{11}E_5^4}{E_1^{14}}+2000q^2
\dfrac{E_2^{11}E_5^{10}}{E_1^{20}}\Bigg).
\end{align*}
Furthermore, for any nonnegative integer $n$,
\begin{align*}f(18750n+15625)&\equiv0~(\textup{mod}~125).
\end{align*}
\end{theorem}

Now we briefly mention a few simple consequences of the generating function \eqref{gen-om-8n4} and our work on this function.

We observe that Waldherr's congruences \eqref{waldherr-cong} can be easily deduced from \eqref{gen-om-8n4} as follows. By the binomial theorem,
\begin{align}\label{binomial}E_1^5\equiv E_5 ~(\textup{mod}~5).\end{align} Therefore, from \eqref{gen-om-8n4}, we have
$$\sum_{n=0}^\infty p_\omega(8n+4)q^n\equiv4\dfrac{E_{10}^{2}E_1}{E_5^2}~(\textup{mod}~5),$$ which can be rewritten with the aid of famous Euler's pentagonal number theorem \cite[p. 12]{spirit}
$$E_1=\sum_{k=-\infty}^\infty(-1)^kq^{k(3k+1)/2},$$as
$$\sum_{n=0}^\infty p_\omega(8n+4)q^n\equiv4\dfrac{E_{10}^{2}}{E_5^2}\sum_{k=-\infty}^\infty(-1)^kq^{k(3k+1)/2}~(\textup{mod}~5).$$
Since $k(3k+1)/2\equiv0,~ 1,~ \textup{or}~2~(\textup{mod}~5)$ only, equating the coefficients of $q^{5n+r}$, $r=3,4$ from both sides of the above, we easily arrive at \eqref{waldherr-cong}.

In this paper, we prove the following representation of the generating function of $p_\omega(40n+12)$.

\begin{theorem}\label{thm-p-omega} We have
\begin{align}\label{gen-om-40n12-main}\sum_{n=0}^\infty p_\omega(40n+12)q^n&=4\Bigg(9\dfrac{E_{2}^{2}E_5}{E_1^{2}}+3975q\dfrac{E_2^3E_{10}^3}{E_1^5}+207425q^2
\dfrac{E_2^4E_{10}^{6}}{E_1^{8}E_5}\notag\\
&\quad+4229000q^3\dfrac{E_{2}^{5}E_{10}^{9}}{E_1^{11}E_5^2}+44850000q^4\dfrac{E_{2}^{6}E_{10}^{12}}{E_1^{14}E_5^3}\notag\\
&\quad+274000000q^5\dfrac{E_{2}^{7}E_{10}^{15}}{E_1^{17}E_5^4}+980000000q^6\dfrac{E_{2}^{8}E_{10}^{18}}{E_1^{20}E_5^5}\notag\\
&\quad+1920000000q^7\dfrac{E_{2}^{9}E_{10}^{21}}{E_1^{23}E_5^6}+1600000000q^8\dfrac{E_{2}^{10}E_{10}^{24}}{E_1^{26}E_5^7}\Bigg).
\end{align}
\end{theorem}

We also deduce the following interesting congruence  recently proved  by Xia \cite{xia-mock}.

\begin{corollary}\label{cor-p-omega}
\noindent For integers $n\geq 0$ and $k\geq0$, we have
\begin{align}\label{p-omega-inf}
p_\omega\left(8\times5^{2k+1}n+\dfrac{7\times5^{2k+1}+1}{3}\right)  \equiv  (-1)^k p_\omega\left(40n+12\right) ~(\textup{mod}~5).
\end{align}
\end{corollary}

The smallest parts function $\textup{spt}(n)$, counting the total number of appearances of the smallest parts in all partitions of $n$, was introduced by Andrews \cite{andrews-crelle}, and the function has received great attention since its introduction. For example, see \cite{andrews-jcta, andrews-spt-rama,  bringmann, dixit, garvan, jennings-jnt1, jennings-jnt, jennings-adv}.

Andrews, Dixit and Yee \cite{andrews-research} also studied the associated smallest parts functions $\textup{spt}_\omega(n)$ and $\textup{spt}_\nu(n)$, which count the number of smallest parts in the partitions enumerated by $p_\omega(n)$ and $p_\nu(n)$, respectively. Of course, $\textup{spt}_\nu(n)=p_\nu(n)$. They proved the congruences
\begin{align}
\label{spt-om-5n3}\textup{spt}_\omega(5n+3)&\equiv  0~(\textup{mod}~5),\\
\label{spt-om-10n7}\textup{spt}_\omega(10n+7)&\equiv  0~(\textup{mod}~5),\\
\label{spt-om-10n9}\textup{spt}_\omega(10n+9)&\equiv  0~(\textup{mod}~5),\\
\label{spt-nu-10n8}\textup{spt}_\nu(10n+8)&=p_\nu(10n+8)\equiv  0~(\textup{mod}~5),
\end{align}
where the first congruence was also independently established by Garvan and Jennings-Shaffer \cite{garvan-jennings}. In fact, Garvan and Jennings-Shaffer \cite{garvan-jennings} introduced a crank-type function that explains the congruence \eqref{spt-om-5n3}. The asymptotic behavior of that crank function was studied by Jang and Kim \cite{jang-kim}. As mentioned earlier, \eqref{spt-nu-10n8} immediately follows from \eqref{gen-nu-10n8}, a fact not possibly noticed by the authors of \cite{andrews-research}. In \cite{andrews-acta}, Andrews, Dixit, Schultz and Yee studied the overpartition analogue of $p_\omega(n)$, namely,  $\overline{p}_\omega(n)$, which counts the number of overpartitions of $n$ such that all odd parts are less than twice the smallest part, and in which the smallest part is always overlined. They also studied $\overline{\textup{spt}}_\omega(n)$, the number of smallest parts in the overpartitions of $n$ in which the smallest part is always overlined and all odd parts are less than twice the smallest part. They found several congruences modulo 2, 3, 4, 5 and 6 for $\overline{p}_\omega(n)$ and $\overline{\textup{spt}}_\omega(n)$. They \cite[Problem 1]{andrews-acta} also raised the question of  relating the generating function of $\overline{p}_\omega(n)$ to modular forms. Recently, the question was answered in affirmative by Bringmann, Jennings-Shaffer and Mahlburg \cite{bringmann-adv}.

Recently, Wang \cite{wang} and Cui, Gu and Hao \cite{cui-gu-hao} also found many new congruences satisfied by $\overline{p}_\omega(n)$, $\textup{spt}_\omega(n)$ and  $\overline{\textup{spt}}_\omega(n)$ modulo powers of 2 and 3. In particular, Wang \cite{wang} derived the following exact generating functions:
\begin{align}
\label{gen-spt-om-2n1}\sum_{n=0}^\infty \textup{spt}_\omega(2n+1)q^n&=\dfrac{E_2^{8}}{E_1^5},\\
\label{gen-overspt-om-2n1}\sum_{n=0}^\infty \overline{\textup{spt}}_\omega(2n+1)q^n&=\dfrac{E_2^{9}}{E_1^6}.
\end{align}

We note that congruences \eqref{spt-om-10n7} and \eqref{spt-om-10n9} can be easily deduced from \eqref{gen-spt-om-2n1} as in the following. Taking congruences modulo 5 in \eqref{gen-spt-om-2n1} and using \eqref{binomial}, we have
\begin{align}\label{cong-w}\sum_{n=0}^\infty \textup{spt}_\omega(2n+1)q^n\equiv\dfrac{E_{10}E_2^{3}}{E_5}~(\textup{mod}~5).
\end{align}Employing Jacobi's identity \cite[p. 14]{spirit}
\begin{align}\label{jacobi}E_1^3=\sum_{k=0}^\infty(-1)^k(2k+1)q^{k(k+1)/2},\end{align} in \eqref{cong-w}, we have
\begin{align}\label{cong-w1}\sum_{n=0}^\infty \textup{spt}_\omega(2n+1)q^n\equiv\dfrac{E_{10}}{E_5}\sum_{k=0}^\infty(-1)^k(2k+1)q^{k(k+1)}~(\textup{mod}~5).
\end{align}
Since $k(k+1)\equiv0,~ 1,~ \textup{or}~2~(\textup{mod}~5)$, equating the coefficients of $q^{5n+r}$, $r=3,4$ from both sides of the above, we easily arrive at \eqref{spt-om-10n7} and \eqref{spt-om-10n9}. Furthermore, we note that $k(k+1)\equiv 1~ (\textup{mod}~5)$ only when $k\equiv 2~ (\textup{mod}~5)$, that is, only when
$2k+1\equiv 0~ (\textup{mod}~5)$. Therefore, equating the coefficients of $q^{5n+1}$ from both sides of the above we arrive at
\begin{align}
\label{cong10n3}\textup{spt}_\omega\left(10n+3\right)\equiv  0~(\textup{mod}~5),
\end{align}
which is, in fact,  contained in \eqref{spt-om-5n3}.

\vspace{.5cm}Wang \cite{wang} offered the following interesting conjecture.

\begin{conjecture}For any integers $k\ge 1$ and $n\ge 0$,
\begin{align}
\textup{spt}_\omega\left(2\cdot 5^{2k-1}n+\dfrac{7\cdot 5^{2k-1}+1}{12}\right)\equiv  0~(\textup{mod}~5^{2k-1}),\\
\textup{spt}_\omega\left(2\cdot 5^{2k}n+\dfrac{11\cdot 5^{2k}+1}{12}\right)\equiv  0~(\textup{mod}~5^{2k}).
\end{align}
\end{conjecture}

The above conjecture has been proved recently with the aid of modular forms by Wang and Yang \cite{wang-yang}.

The cases $k=1$ and $2$ of the above congruences are \eqref{cong10n3},
\begin{align}
\label{cong50n23}\textup{spt}_\omega\left(50n+23\right)\equiv  0~(\textup{mod}~5^{2}),\\
\label{cong250n73}\textup{spt}_\omega\left(250n+73\right)\equiv  0~(\textup{mod}~5^{3}),\\\intertext{and}
\label{cong1250n573}\textup{spt}_\omega\left(1250n+573\right)\equiv  0~(\textup{mod}~5^{4}).
\end{align}
In this paper, we find the following exact generating functions of $\textup{spt}_\omega\left(10n+3\right)$ and $\textup{spt}_\omega\left(50n+23\right)$.
\begin{theorem}\label{thm-spt-main} We have
\begin{align}\label{gen-spt-10n3-main}&\sum_{n=0}^\infty\textup{spt}_\omega\left(10n+3\right)q^n=5\left(E_1^2E_5+6q\dfrac{E_2E_{10}^3}{E_1}
+25q\dfrac{E_2^8E_5^5}{E_1^{10}}\right),\\
             &\sum_{n=0}^\infty\textup{spt}_\omega\left(50n+23\right)q^n\notag\\
&=25\Bigg(E_1E_5^2+6q\dfrac{E_2E_5E_{10}^3}{E_1^2}+25\left(2q\dfrac{E_1^3E_{10}^4}{E_2^4}+25q^3\dfrac{E_1^3E_{10}^{10}}{E_2^{10}}\right)\notag\\
&\quad+50 \Bigg(\dfrac{E_1^{16}E_{10}^{3}}{E_2^{15}E_5}+150q
\dfrac{E_1^{13}E_{10}^{6}}{E_2^{14}E_5^2}+5650q^2\dfrac{E_1^{10}E_{10}^{9}}{E_2^{13}E_5^3}+101825q^3\dfrac{E_1^{7}E_{10}^{12}}{E_2^{12}E_5^4}
\notag\\
&\quad+1068125q^4\dfrac{E_1^{4}E_{10}^{15}}{E_2^{11}E_5^5}+7042500q^5\dfrac{E_1E_{10}^{18}}{E_2^{10}E_5^6}+29800000q^6\dfrac{E_{10}^{21}}{E_1^{2}E_2^{9}E_5^7}
\notag\\
&\quad+79000000q^7\dfrac{E_{10}^{24}}{E_1^{5}E_2^{8}E_5^8}+120000000q^8\dfrac{E_{10}^{27}}{E_1^{8}E_2^{7}E_5^9}
+80000000q^9\dfrac{E_{10}^{30}}{E_1^{11}E_2^{6}E_5^{10}}\Bigg)\notag\\\label{gen-spt-50n23-main}
&\quad+ 625q \Bigg(63\dfrac{E_2^8E_5^6}{E_1^{11}}+6500q\dfrac{E_2^8E_5^{12}}{E_1^{17}}+196875q^2\dfrac{E_2^8E_5^{18}}{E_1^{23}}+2343750q^3\dfrac{E_2^8E_5^{24}}{E_1^{29}}\notag\\
&\quad+9765625q^4\dfrac{E_2^8E_5^{30}}{E_1^{35}}\Bigg)\Bigg).
\end{align}
\end{theorem}
Note that the congruences \eqref{cong10n3} and \eqref{cong50n23} follow trivially from the above theorem. In this paper, we also deduce \eqref{cong250n73} and \eqref{cong1250n573}.

There are several congruences for $\overline{\textup{spt}}_\omega(n)$ modulo 11 and powers of 2 and 3 (For example, see \cite{andrews-acta, wang, cui-gu-hao}). But for modulo 5, to our knowledge, the following congruence, found by Andrews et al. \cite{andrews-acta}, is the only available one:
$$\overline{\textup{spt}}_\omega\left(10n+6\right)\equiv  0~(\textup{mod}~5).$$
In this paper, we present the following exact generating function of $\overline{\textup{spt}}_\omega\left(10n+5\right)$ and $\overline{\textup{spt}}_\omega\left(50n+25\right)$.
\begin{theorem}\label{thm-sptbar}We have
\begin{align}\label{gen-over-spt-10n5-main}\sum_{n=0}^\infty \overline{\textup{spt}}_\omega\left(10n+5\right)q^n&=18E_1^2E_{10}+720q\dfrac{E_2E_{10}^4}{E_1E_5}+7625q^2
\dfrac{E_2^2E_{10}^{7}}{E_1^{4}E_5^2}+32500q^3\dfrac{E_{2}^{3}E_{10}^{10}}{E_1^{7}E_5^3}\notag\\
&\quad+50000q^4\dfrac{E_{2}^{4}E_{10}^{13}}{E_1^{10}E_5^4}\\\intertext{and}
\label{gen-over-spt-50n25-main}\sum_{n=0}^\infty \overline{\textup{spt}}_\omega\left(50n+25\right)q^n
&=8327 E_2E_5^{2}
+28312350q\dfrac{E_2^{2}E_5E_{10}^{3}}{E_1^{3}}
+7557865625q^2\dfrac{E_2^{3}E_{10}^{6}}{E_1^{6}}
\notag\\
&\quad+678027312500q^3\dfrac{E_2^{4}E_{10}^{9}}{E_1^{9}E_5}+3072484775\times10^4q^4\dfrac{E_2^{5}E_{10}^{12}}{E_1^{12}E_5^{2}}\notag\\
&\quad+84314744\times10^7q^5\dfrac{E_2^{6}E_{10}^{15}}{E_1^{15}E_5^{3}}
+154486601\times10^8q^6\dfrac{E_2^{7}E_{10}^{18}}{E_1^{18}E_5^{4}}\notag\\
&\quad+20019467\times10^{10}q^7\dfrac{E_2^{8}E_{10}^{21}}{E_1^{21}E_5^{5}}
+18998414\times10^{11}q^8\dfrac{E_2^{9}E_{10}^{24}}{E_1^{24}E_5^{6}}
\notag\\
&\quad+134691824\times10^{11}q^9\dfrac{E_2^{10}E_{10}^{27}}{E_1^{27}E_5^{7}}+71952464\times10^{12}q^{10}\dfrac{E_2^{11}E_{10}^{30}}{E_1^{30}E_5^{8}}\notag\\
&\quad+28911776\times10^{13}q^{11}\dfrac{E_2^{12}E_{10}^{33}}{E_1^{33}E_5^{9}}+862432\times10^{15}q^{12}\dfrac{E_2^{13}E_{10}^{36}}{E_1^{36}E_5^{10}}\notag\\
&\quad+1856\times10^{18}q^{13}\dfrac{E_2^{14}E_{10}^{39}}{E_1^{39}E_5^{11}}+272896\times10^{16}q^{14}\dfrac{E_2^{15}E_{10}^{42}}{E_1^{42}E_5^{12}}\notag\\
&\quad
+24576\times10^{17}q^{15}\dfrac{E_2^{16}E_{10}^{45}}{E_1^{45}E_5^{13}}+1024\times10^{18}q^{16}\dfrac{E_2^{17}E_{10}^{48}}{E_1^{48}E_5^{14}}.
\end{align}
\end{theorem}

As a consequence, we also deduce the following new congruences.
\begin{corollary}\label{cor-sptbar}Let $\ell \in\{1,2,3,4,5,6\}$. Then for integers $k\geq0$ and $n\geq0$, we have

\begin{align}\label{spt-gencong}&\overline{\textup{spt}}_\omega\left(5^{2k+\ell-1}\left(10n+5\right)\right)\equiv \overline{\textup{spt}}_\omega\left(5^{\ell-1}\left(10n+5\right)\right)~(\textup{mod}~5^{\ell}).
\end{align}

Furthermore, we have
\begin{align}\label{spt-cor2-mod5}&\overline{\textup{spt}}_\omega\left(5^{2}\left(10n+3\right)\right)\equiv \overline{\textup{spt}}_\omega\left(5^{2}\left(10n+7\right)\right)\equiv0~(\textup{mod}~5^2),\\
\label{spt-cor2-mod125}&\overline{\textup{spt}}_\omega\left(5^{4}\left(10n+3\right)\right)\equiv \overline{\textup{spt}}_\omega\left(5^{4}\left(10n+7\right)\right)\equiv0~(\textup{mod}~5^4),\\
\label{spt-cor2-mod3125}&\overline{\textup{spt}}_\omega\left(5^{6}\left(10n+3\right)\right)\equiv \overline{\textup{spt}}_\omega\left(5^{6}\left(10n+7\right)\right)\equiv0~(\textup{mod}~5^6).
\end{align}
\end{corollary}

It is worthwhile to mention that congruences similar to the above for modulo $5^\ell$ for some $\ell>6$ might be deduced similarly from \eqref{gen-over-spt-50n25-main}. We propose the following general conjecture.

\begin{conjecture} For integers $\ell\ge 1$, $k\geq0$ and $n\geq0$, we have
\begin{align*}&\overline{\textup{spt}}_\omega\left(5^{2k+\ell-1}\left(10n+5\right)\right)\equiv \overline{\textup{spt}}_\omega\left(5^{\ell-1}\left(10n+5\right)\right)~(\textup{mod}~5^{\ell})\\\intertext{and}
&\overline{\textup{spt}}_\omega\left(5^{2\ell}\left(10n+3\right)\right)\equiv \overline{\textup{spt}}_\omega\left(5^{2\ell}\left(10n+7\right)\right)\equiv0~(\textup{mod}~5^{2\ell}).
\end{align*}
\end{conjecture}

We employ some well-known identities for the Rogers-Ramanujan continued fraction, which is defined by
 \begin{equation}\label{rq}
R(q):=\dfrac{q^{1/5}}{1}\+\dfrac{q}{1}\+\dfrac{q^2}{1}\+\dfrac{q^3}{1}\+\cds=q^{1/5}\dfrac{(q;q^5)_\infty(q^4;q^5)_\infty}{(q^2;q^5)_\infty(q^3;q^5)_\infty},~|q|<1,
\end{equation}
and some identities developed in our paper \cite{baruah-begum}.

We organize the paper in the following way. In the next section, we present some useful definitions and lemmas that will be used in the subsequent sections. In Section \ref {sec3}, we prove Theorem \ref{thm-nu-50n8}, and deduce \eqref{p-nu-250n208} and Corollary \ref{cor-nu-50n8}. In Section \ref{sec4}, we prove Theorem \ref{thm-p-omega} and Corollary \ref{cor-p-omega}.
In Section \ref{sec5}, we prove Theorem \ref{thm-spt-main} and the congruences \eqref{cong250n73} and \eqref{cong1250n573}. In the last section, we prove Theorem \ref{thm-sptbar} and Corollary \ref{cor-sptbar}.
\section{\textbf{Some useful definitions and lemmas}}
Ramanujan's theta functions $\varphi(-q)$ and $\psi(q)$ are defined, for $|q|<1$, by
\begin{align}\label{varphi}
\varphi(-q)&:=\sum_{j=-\infty}^\infty (-1)^j
q^{j^2}=(q; q^2)^2_\infty(q^2; q^2)_\infty=\dfrac{E_1^2}
{E_2}\\\intertext{and}
\label{psi}\psi(q)&:=\sum_{j=0}^\infty q^{j(j+1)/2}=\dfrac{(q^2; q^2)_\infty}{(q; q^2)_\infty}=\dfrac{E_2^2}
{E_1},
\end{align}
where the product representations arise from Jacobi's famous triple product identity \cite[p. 35, Entry 19]{bcb3}.

In the following lemma we record some well-known identities (See Berndt's books \cite[p. 40]{bcb3} and \cite[p. 165]{spirit}), where the first three are 5-dissections of $E_1$, $1/E_1$, and $\varphi(-q)$, respectively.
\begin{lemma}\label{dissect} If ~ $T(q):=\dfrac{q^{1/5}}{R(q)}=\dfrac{(q^2;q^5)_\infty(q^3;q^5)_\infty}{(q;q^5)_\infty(q^4;q^5)_\infty}$, then
\begin{align}
\label{E1}E_1&=E_{25}\left(T(q^5)-q-\dfrac{q^{2}}{T(q^5)}\right),\\
\label{1byE1}\dfrac{1}{E_1}&= \dfrac{E_{25}^5}{E_5^6}\Bigg(T(q^5)^4+qT(q^5)^3+2q^2T(q^5)^2+3q^3T(q^5)+5q^4-\dfrac{3q^5}{T(q^5)}
\notag\\
&\quad+\dfrac{2q^6}{T(q^5)^2}-\dfrac{q^7}{T(q^5)^3}+\dfrac{q^8}{T(q^5)^4}\Bigg),\\
\label{phi-5dissect}\varphi(-q)&=\dfrac{E_1^2}
{E_2}=\dfrac{E_{25}^2}
{E_{50}}-2q(q^{15},q^{35},q^{50};q^{50})_\infty +2q^4 (q^{5},q^{45},q^{50};q^{50})_\infty,\\\intertext{and}
\label{E1-6}11q&+\dfrac{E_1^6}{E_5^6}= T(q)^5-\dfrac{q^2}{T(q)^5}.
\end{align}
\end{lemma}

In the following two lemmas, we recall some useful results from our recent paper \cite{baruah-begum}.
\begin{lemma}\label{lemmax-y}\textup{(Baruah and Begum \cite[Lemma 1.3]{baruah-begum})}\label{xy} If $x=T(q)$ and $y=T(q^2)$, then
\begin{align}
\label{xy2}xy^2-\dfrac{q^2}{xy^2}&=K,\\
\label{x2-by-y}\dfrac{x^2}{y}-\dfrac{y}{x^2}&=\dfrac{4q}{K},\\
\label{y3-by-x}\dfrac{y^3}{x}+q^2\dfrac{x}{y^3}&=K+\dfrac{4q^2}{K}-2q,\\
\label{x3y}x^3y+\dfrac{q^2}{x^3y}&=K+\dfrac{4q^2}{K}+2q,
\end{align}
where $K=(E_2E_5^5)/(E_1E_{10}^5).$
\end{lemma}
\begin{lemma}\label{AB}\textup{(Baruah and Begum \cite[Eqs. (2.6), (2.7), (2.29)]{baruah-begum})}
\begin{align}\label{A-4qB}
\dfrac{E^5_5}{E_1^{4}E^{3}_{10}}&=\dfrac{E_5}{E_2^{2}E_{10}}+4q\dfrac{E^{2}_{10}}{E_1^{3}E_2},\\
\label{A-qB}
\dfrac{E^3_2E_5^2}{E_1^{2}E^2_{10}}&=\dfrac{E^5_5}{E_1E^{3}_{10}}
+q\dfrac{E^{2}_{10}}{E_2},\\
\label{A-5qB}
\dfrac{E^3_2E_5^2}{E_1^5E^{2}_{10}}&=\dfrac{E_5}{E_2^2E_{10}}+5q\dfrac{E^2_{10}}{E^3_1E_2}.
\end{align}
\end{lemma}

\section{\textbf{Proofs of Theorem \ref{thm-nu-50n8}, Eq. \eqref{p-nu-250n208} and Corollary \ref{cor-nu-50n8}}}\label{sec3}

\begin{pf-thm-nu-50n8}
Employing \eqref{A-4qB} successively in \eqref{gen-nu-10n8},  we see that
\begin{align*}\sum_{n=0}^\infty p_\nu(10n+8)q^n&=5\left(\dfrac{E_{10}^{3}}{E_5^{2}}+4q\dfrac{E_2E_{10}^6}{E_1^{3}E_5^3}\right)\\
&=5\left(\dfrac{E_{10}^{3}}{E_5^{2}}+4q\dfrac{E_1E_{10}^8}{E_2E_5^7}+16q^2\dfrac{E_{10}^{11}}{E_1^{2}E_5^8}\right)
\end{align*}
Employing \eqref{E1} and \eqref{1byE1} in the above, extracting the terms involving $q^{5n}$, and then replacing $q^5$ by $q$, we find that
\begin{align*}\sum_{n=0}^\infty p_\nu(50n+8)q^n&=5\dfrac{E_{2}^{3}}{E_1^{2}}+20q\dfrac{E_2^2E_5E_{10}^5}{E_1^7}\Bigg(2\left(xy^2-\dfrac{q^2}{xy^2}\right)
-\left(\dfrac{y^3}{x}+q^2\dfrac{x}{y^3}\right)-5q\Bigg)\\
&\quad+400q\dfrac{E_{2}^{11}E_5^{10}}{E_1^{20}}\left(2\left(x^5-\dfrac{q^2}{x^5}\right)+3q\right).
\end{align*}
Employing \eqref{E1-6}, \eqref{xy2} and \eqref{y3-by-x} in the above, we find that
\begin{align*}\sum_{n=0}^\infty p_\nu(50n+8)q^n&=5\dfrac{E_{2}^{3}}{E_1^{2}}+20q\dfrac{E_2^3E_5^6}{E_1^8}-60q^2\dfrac{E_2^2E_5E_{10}^5}{E_1^7}-80q^3
\dfrac{E_2E_{10}^{10}}{E_1^6E_5^4}\notag\\
&\quad+800q\dfrac{E_{2}^{11}E_5^{4}}{E_1^{14}}+10000q^2\dfrac{E_{2}^{11}E_5^{10}}{E_1^{20}}.
\end{align*}
We reduce the above into the form \eqref{gen-nu-50n8-main} with the aid of \eqref{A-4qB} and \eqref{A-qB} as shown below:
\begin{align*}\sum_{n=0}^\infty p_\nu(50n+8)q^n&=5\dfrac{E_{2}^{3}}{E_1^{2}}+20q\dfrac{E_2^3E_5E_{10}^3}{E_1^7}\left(\dfrac{E^5_5}{E_1E^{3}_{10}}
+q\dfrac{E^{2}_{10}}{E_2}\right)-80q^2\dfrac{E_2^2E_{10}^8}{E_1^6E_5^4}\\
&\quad\times\left(\dfrac{E^5_5}{E_1E^{3}_{10}}
+q\dfrac{E^{2}_{10}}{E_2}\right)+800q\dfrac{E_{2}^{11}E_5^{4}}{E_1^{14}}+10000q^2\dfrac{E_{2}^{11}E_5^{10}}{E_1^{20}}\\
&=5\dfrac{E_{2}^{3}}{E_1^{2}}+20q\dfrac{E_2^6E_5^3E_{10}}{E_1^9}-80q^2\dfrac{E_2^5E_{10}^6}{E_1^8E_5^2}+800q\dfrac{E_{2}^{11}E_5^{4}}{E_1^{14}}\\
&\quad+10000q^2\dfrac{E_{2}^{11}E_5^{10}}{E_1^{20}}\\
&=5\dfrac{E_{2}^{3}}{E_1^{2}}+20q\dfrac{E_2^6E_{10}^4}{E_1^5E_5^2}\left(\dfrac{E^5_5}{E_1^4E_{10}^{3}}-4q\dfrac{E_{10}^2}{E_1^3E_2}\right)
+800q\dfrac{E_{2}^{11}E_5^{4}}{E_1^{14}}\\
&\quad+10000q^2\dfrac{E_{2}^{11}E_5^{10}}{E_1^{20}}\\
&=5\dfrac{E_{2}^{3}}{E_1^{2}}+20q\dfrac{E_2^4E_{10}^3}{E_1^5E_5}
+800q\dfrac{E_{2}^{11}E_5^{4}}{E_1^{14}}
+10000q^2\dfrac{E_{2}^{11}E_5^{10}}{E_1^{20}}\\
&=5\dfrac{E_{2}^{5}E_5^4}{E_1^6E_{10}^2}
+800q\dfrac{E_{2}^{11}E_5^{4}}{E_1^{14}}
+10000q^2\dfrac{E_{2}^{11}E_5^{10}}{E_1^{20}}.
\end{align*}
\end{pf-thm-nu-50n8}

\begin{pf-cor1-nu-50n8}
Employing \eqref{binomial} in \eqref{gen-nu-50n8-main}, then using \eqref{p-n}, we have
\begin{align*}\sum_{n=0}^\infty p_\nu(50n+8)q^n&\equiv5~\dfrac{E_5^3}{E_1E_{10}}\equiv5~\dfrac{E_5^3}{E_{10}}~\sum_{k=0}^\infty p(k)q^k~(\textup{mod}~25).
\end{align*}
Extracting the coefficients of $q^{5n+4}$, and then using \eqref{rama-beautiful}, we easily arrive at \eqref{p-nu-250n208}.
\end{pf-cor1-nu-50n8}

\begin{pf-cor2-nu-50n8}
From \eqref{gen-nu-50n8-main}, we have
\begin{align*}\sum_{n=0}^\infty p_\nu(50n+8)q^n&\equiv5\Bigg(\dfrac{E_{2}^{5}E_5^4}{E_1^{6}E_{10}^2}+160q\dfrac{E_2^{11}E_5^4}{E_1^{14}}\Bigg)~(\textup{mod}~125).
\end{align*}
Employing \eqref{A-5qB} in the above, we see that
\begin{align}\label{gen-nu-1250n1280-2}\sum_{n=0}^\infty p_\nu(50n+8)q^n&\equiv5\Bigg(\dfrac{E_5^3}{E_1E_{10}}+165q\dfrac{E_2E_5^2E_{10}^2}{E_1^4}\Bigg)\notag\\
&\equiv5\Bigg(\dfrac{E_5^3}{E_1E_{10}}+165qE_1E_2E_5E_{10}^2\Bigg)~(\textup{mod}~125).
\end{align}
Employing \eqref{E1} and \eqref{1byE1} in the above, extracting the terms involving $q^{5n+4}$ from both sides, dividing by $q^4$, then replacing $q^5$ by $q$, and then using \eqref{rama-beautiful}, we obtain
\begin{align}\label{gen-nu-1250n1280-22}\sum_{n=0}^\infty p_\nu(250n+208)q^n&\equiv25\Bigg(\dfrac{E_5^5}{E_1^3E_{2}}+33E_1E_2^2E_5E_{10}\Bigg)\notag\\
&\equiv25\Bigg(E_5^4\cdot\dfrac{E_1^2}{E_2}+33E_1E_2^2E_5E_{10}\Bigg)~(\textup{mod}~125).
\end{align}
Employing \eqref{E1} and \eqref{phi-5dissect} in \eqref{gen-nu-1250n1280-22}, extracting the terms involving $q^{5n}$ from both sides, and then replacing $q^5$ by $q$, we find that
\begin{align}\label{gen-nu-1250n1280-3}\sum_{n=0}^\infty p_\nu(1250n+208)q^n&\equiv25\left(\dfrac{E_1^4E_5^2}{E_{10}}+33E_1E_2E_5E_{10}^2\left(xy^2+q-\dfrac{q^2}{xy^2}\right)\right)\notag\\
&\equiv25\left(\dfrac{E_5^3}{E_1E_{10}}+33E_1E_2E_5E_{10}^2\left(xy^2+q-\dfrac{q^2}{xy^2}\right)\right)~(\textup{mod}~125),
\end{align}
where $x$ and $y$ are as defined in Lemma \ref{xy}.

Now, from \cite[p. 56]{ramalost} (\cite[p. 35, Entry 1.8.2]{lost1}), we note that if $k=qR(q)R^2(q^2)$ and $k\leq\sqrt{5}-2$, then
\begin{equation*}
\dfrac{\psi^2(q)}{q\psi^2(q^5)}=\dfrac{1+k-k^2}{k},
\end{equation*}
which can be seen to be equivalent to
\begin{equation}\label{k2}
xy^2+q-\dfrac{q^2}{xy^2}=\dfrac{\psi^2(q)}{\psi^2(q^5)}=\dfrac{E_2^4E_5^2}{E_1^2E_{10}^4},
\end{equation}
where the last equality is by \eqref{psi}.

Using \eqref{k2} in \eqref{gen-nu-1250n1280-3} and then applying \eqref{binomial}, we obtain
\begin{align*}\sum_{n=0}^\infty p_\nu(1250n+208)q^n&\equiv25\times 34~\dfrac{E_5^3}{E_1E_{10}}\\
&\equiv25\times 34~\dfrac{E_5^3}{E_{10}}\sum_{n=0}^\infty p(n)q^n~(\textup{mod}~125).
\end{align*}
Equating the coefficients of $q^{5n+4}$ from both sides of the above, and then applying \eqref{rama-beautiful}, we arrive at
\begin{align*}p_\nu(6250n+5208)&\equiv0~(\textup{mod}~125),
\end{align*}
to finish the proof.
\end{pf-cor2-nu-50n8}

\section{\textbf{Proofs of Theorem \ref{thm-p-omega} and Corollary \ref{cor-p-omega}}}\label{sec4}
\begin{pf-thm-p-omega}
Employing \eqref{A-5qB}
in \eqref{gen-om-8n4} successively, we have
\begin{align*}\sum_{n=0}^\infty p_\omega(8n+4)q^n&=4\Bigg(\dfrac{E_2^5E_{10}}{E_1^4E_5}+5q\dfrac{E_{2}^{6}E_{10}^{4}}{E_1^{7}E_5^{2}}\Bigg)\\
&=4\Bigg(\dfrac{E_1E_{10}^{2}}{E_5^{2}}+5q\dfrac{E_2E_{10}^5}{E_1^2E_5^3}+5q\left(\dfrac{E_2E_{10}^5}{E_1^2E_5^3}+5q\dfrac{E_2^2E_{10}^8}{E_1^5E_5^4}\right)\Bigg)\\
&=4\Bigg(\dfrac{E_1E_{10}^{2}}{E_5^{2}}+10q\dfrac{E_2E_{10}^5}{E_1^2E_5^3}+25q^2\dfrac{E_{2}^{2}E_{10}^{8}}{E_1^{5}E_5^{4}}\Bigg).
\end{align*}
Now applying \eqref{A-4qB}  in the above successively, we see that
\begin{align}\label{gen-om-40n12-1}\sum_{n=0}^\infty p_\omega(8n+4)q^n&=4\Bigg(\dfrac{E_1E_{10}^{2}}{E_5^{2}}+10q\left(\dfrac{E_1^2E_{10}^7}{E_2E_5^7}+4q\dfrac{E_{10}^{10}}{E_1E_5^{8}}\right)\notag\\
&\quad+25q^2\left(\dfrac{E_{10}^{10}}{E_1E_5^{8}}+4q\dfrac{E_2E_{10}^{13}}{E_1^4E_5^{9}}\right)\Bigg)\notag\\
&=4\Bigg(\dfrac{E_1E_{10}^{2}}{E_5^{2}}+10q\dfrac{E_1^2E_{10}^7}{E_2E_5^7}+65q^2\dfrac{E_{10}^{10}}{E_1E_5^{8}}
+100q^3\dfrac{E_{10}^{15}}{E_2E_5^{13}}\notag\\
&\quad+400q^4\dfrac{E_{10}^{18}}{E_1^3E_5^{14}}\Bigg).
\end{align}
Employing \eqref{E1} and \eqref{1byE1} in the above, extracting the terms involving $q^{5n+1}$, dividing both sides by $q$, and then replacing $q^5$ by $q$, we obtain
\begin{align*}\sum_{n=0}^\infty p_\omega(40n+12)q^n&=-4\dfrac{E_{2}^{2}E_5}{E_1^{2}}+40\dfrac{E_2E_5^2E_{10}^5}{E_1^7}\left(x^2y^4+\dfrac{q^4}{x^2y^4}\right)+q\Bigg(1300\dfrac{E_{2}^{10}E_5^{5}}{E_1^{14}}\notag\\
&\quad+14400\dfrac{E_{2}^{18}E_5^{15}}{E_1^{32}}\left(x^{10}+\dfrac{q^4}{x^{10}}\right)
-160\dfrac{E_{2}E_5^{2}E_{10}^5}{E_{1}^7}\left(xy^2-\dfrac{q^2}{xy^2}\right)\\
&\quad+80\dfrac{E_{2}E_5^{2}E_{10}^5}{E_{1}^7}\left(\dfrac{y^3}{x}+q^2\dfrac{x}{y^3}\right)\Bigg) +q^2\Bigg(2000\dfrac{E_{2}^9E_{10}^5}{E_{1}^{13}}\notag\\
&\quad-200\dfrac{E_{2}E_5^{2}E_{10}^5}{E_{1}^7}+283200\dfrac{E_{2}^{18}E_5^{15}}{E_1^{32}}\left(x^{5}-\dfrac{q^2}{x^{5}}\right)\\
&\quad-120\dfrac{E_{2}E_5^{2}E_{10}^5}{E_{1}^7}\left(\dfrac{x^2}{y}-\dfrac{y}{x^2}\right)\Bigg)
+113600q^3\dfrac{E_{2}^{18}E_5^{15}}{E_1^{32}}.
\end{align*}
With the aid of \eqref{xy2} -- \eqref{y3-by-x}, the above reduces to
\begin{align}\label{omega-main}&\sum_{n=0}^\infty p_\omega(40n+12)q^n\notag\\
&=-4\dfrac{E_{2}^{2}E_5}{E_1^{2}}+40\dfrac{E_2^3E_5^{12}}{E_1^9E_{10}^5}-80q\dfrac{E_2^2E_5^7}{E_1^8}
-280q^2\dfrac{E_2E_5^2E_{10}^5}{E_1^7}-160q^3\dfrac{E_{10}^{10}}{E_1^6E_5^3}\notag\\
&\quad+1300q\dfrac{E_{2}^{10}E_5^{5}}{E_1^{14}}+2000q^2\dfrac{E_{2}^{9}E_{10}^{5}}{E_1^{13}}
+14400q\dfrac{E_{2}^{20}E_5^{25}}{E_{1}^{34}E_{10}^{10}}+398400q^2\dfrac{E_{2}^{19}E_5^{20}}{E_{1}^{33}E_{10}^5}\notag\\
&\quad+1736000q^3 \dfrac{E_{2}^{18}E_{5}^{15}}{E_{1}^{32}}+3648000q^4\dfrac{E_{2}^{17}E_5^{10}E_{10}^5}{E_{1}^{31}}+7296000q^5\dfrac{E_{2}^{16}E_5^{5}E_{10}^{10}}{E_1^{30}}\notag\\
&\quad+3686400q^6\dfrac{E_{2}^{15}E_{10}^{15}}{E_{1}^{29}}+3686400q^7\dfrac{E_{2}^{14}E_{10}^{20}}{E_1^{28}
E_5^5}\notag\\
&=-4\dfrac{E_{2}^{2}E_5}{E_1^{2}}+40\dfrac{E_2^2E_{10}}{E_1^4E_{5}^3}\left(\dfrac{E^5_5}{E_1^{4}E^{3}_{10}}-4q\dfrac{E^{2}_{10}}{E_1^{3}E_2}\right)
\left(\dfrac{E^5_5}{E_1E_{10}^3}+q\dfrac{E_{10}^2}{E_2}\right)\notag\\
&\quad+A(q)+B(q)\notag\\
&=-4\dfrac{E_{2}^{2}E_5}{E_1^{2}}+40\dfrac{E_2^7E_5^2}{E_1^7E_{10}}+A(q)+B(q)\notag\\
&=36\dfrac{E_{2}^{2}E_5}{E_1^{2}}+200q\dfrac{E_2^3E_{10}^3}{E_1^5}+A(q)+B(q),
\end{align}
where we also applied \eqref{A-4qB} -- \eqref{A-5qB}, and, for convenience, set
\begin{align*}A(q)&:=1300q\dfrac{E_{2}^{10}E_5^{5}}{E_1^{14}}+2000q^2\dfrac{E_{2}^{9}E_{10}^{5}}{E_1^{13}}\\\intertext{and}
B(q)&:=14400q\dfrac{E_{2}^{20}E_5^{25}}{E_{1}^{34}E_{10}^{10}}+398400q^2\dfrac{E_{2}^{19}E_5^{20}}{E_{1}^{33}E_{10}^5} +1736000q^3 \dfrac{E_{2}^{18}E_{5}^{15}}{E_{1}^{32}}\notag\\
&\quad+3648000q^4\dfrac{E_{2}^{17}E_5^{10}E_{10}^5}{E_{1}^{31}}+7296000q^5\dfrac{E_{2}^{16}E_5^{5}E_{10}^{10}}{E_1^{30}}
+3686400q^6\dfrac{E_{2}^{15}E_{10}^{15}}{E_{1}^{29}}\\
&\quad+3686400q^7\dfrac{E_{2}^{14}E_{10}^{20}}{E_1^{28}
E_5^5}.\end{align*}

Now, by \eqref{A-4qB}, we have
\begin{align}\label{Aq}A(q)&=1300q\dfrac{E_2^{10}E_{10}^3}{E_1^{10}}\Bigg(\dfrac{E^5_5}{E_1^{4}E^{3}_{10}}-4q\dfrac{E^{2}_{10}}{E_1^{3}E_2}\Bigg)
+7200q^2\dfrac{E_{2}^{9}E_{10}^{5}}{E_1^{13}}\notag\\
&=1300q\dfrac{E_2^{8}E_5E_{10}^{2}}{E_1^{10}}
+7200q^2\dfrac{E_2^{9}E_{10}^{5}}{E_1^{13}}\intertext{and}
B(q)&=14400q\dfrac{E_2^{20}E_5^{20}}{E_1^{30}E_{10}^7}\Bigg(\dfrac{E^5_5}{E_1^{4}E^{3}_{10}}-4q\dfrac{E^{2}_{10}}{E_1^{3}E_2}\Bigg)
+456000q^2\dfrac{E_2^{19}E_5^{15}}{E_1^{29}E_{10}^2}\Bigg(\dfrac{E^5_5}{E_1^{4}E^{3}_{10}}-4q\dfrac{E^{2}_{10}}{E_1^{3}E_2}\Bigg)\notag\\
&\quad+3560000q^3\dfrac{E_2^{18}E_5^{10}E_{10}^3}{E_1^{28}}\Bigg(\dfrac{E^5_5}{E_1^{4}E^{3}_{10}}-4q\dfrac{E^{2}_{10}}{E_1^{3}E_2}\Bigg)
+17888000q^4\dfrac{E_2^{17}E_5^{5}E_{10}^8}{E_1^{27}}\notag\\
&\quad\times\Bigg(\dfrac{E^5_5}{E_1^{4}E^{3}_{10}}-4q\dfrac{E^{2}_{10}}{E_1^{3}E_2}\Bigg)+78848000q^5\dfrac{E_2^{16}E_{10}^{13}}{E_1^{26}}
\Bigg(\dfrac{E^5_5}{E_1^{4}E^{3}_{10}}-4q\dfrac{E^{2}_{10}}{E_1^{3}E_2}\Bigg)\notag\\
&\quad+319078400q^6\dfrac{E_2^{15}E_{10}^{18}}{E_1^{25}E_5^{5}}\Bigg(\dfrac{E^5_5}{E_1^{4}E^{3}_{10}}-4q\dfrac{E^{2}_{10}}{E_1^{3}E_2}\Bigg)
+1280000000q^7\dfrac{E_2^{14}E_{10}^{20}}{E_1^{28}E_5^5}\notag\\
&=14400q\dfrac{E_2^{18}E_5^{21}}{E_1^{30}E_{10}^8}+456000q^2\dfrac{E_2^{17}E_5^{16}}{E_1^{29}E_{10}^3}
+3560000q^3\dfrac{E_2^{16}E_5^{11}E_{10}^2}{E_1^{28}}\notag\\
&\quad+17888000q^4\dfrac{E_2^{15}E_5^{6}E_{10}^7}{E_1^{27}}+78848000q^5\dfrac{E_2^{14}E_5E_{10}^{12}}{E_1^{26}}
+319078400q^6\dfrac{E_2^{13}E_{10}^{17}}{E_1^{25}E_5^{4}}\notag\\
&\quad+1280000000q^7\dfrac{E_2^{14}E_{10}^{20}}{E_1^{28}E_5^5}\notag\end{align}
Applying \eqref{A-4qB} five more steps as in the above for $B(q)$, we obtain
\begin{align}\label{Bq}
B(q)&=14400q\dfrac{E_2^{8}E_5E_{10}^{2}}{E_1^{10}}
+744000q^2\dfrac{E_2^{9}E_{10}^{5}}{E_1^{13}}
+13160000q^3\dfrac{E_2^{10}E_{10}^{8}}{E_1^{16}E_5}\notag\\
&\quad+113600000q^4\dfrac{E_2^{11}E_{10}^{11}}{E_1^{19}E_5^{2}}+528000000q^5\dfrac{E_2^{12}E_{10}^{14}}{E_1^{22}E_5^{3}}
+1280000000q^6\dfrac{E_2^{13}E_{10}^{17}}{E_1^{25}E_5^{4}}\notag\\
&\quad+1280000000q^7\dfrac{E_2^{14}E_{10}^{20}}{E_1^{28}E_5^{5}}.
\end{align}
Therefore,
\begin{align*}
A(q)+B(q)&=15700q\dfrac{E_2^{8}E_5E_{10}^{2}}{E_1^{10}}
+751200q^2\dfrac{E_2^{9}E_{10}^{5}}{E_1^{13}}
+13160000q^3\dfrac{E_2^{10}E_{10}^{8}}{E_1^{16}E_5}\\
&\quad+113600000q^4\dfrac{E_2^{11}E_{10}^{11}}{E_1^{19}E_5^{2}}+528000000q^5\dfrac{E_2^{12}E_{10}^{14}}{E_1^{22}E_5^{3}}
\\
&\quad+1280000000q^6\dfrac{E_2^{13}E_{10}^{17}}{E_1^{25}E_5^{4}}+1280000000q^7\dfrac{E_2^{14}E_{10}^{20}}{E_1^{28}E_5^{5}}.
\end{align*}
Employing \eqref{A-5qB} in the above, we obtain
\begin{align*}
A(q)+B(q)&=15700q\dfrac{E_2^3E_{10}^3}{E_1^5}+829700q^2
\dfrac{E_2^4E_{10}^{6}}{E_1^{8}E_5}+16916000q^3\dfrac{E_{2}^{5}E_{10}^{9}}{E_1^{11}E_5^2}\notag\\
&\quad+179400000q^4\dfrac{E_{2}^{6}E_{10}^{12}}{E_1^{14}E_5^3}+1096000000q^5\dfrac{E_{2}^{7}E_{10}^{15}}{E_1^{17}E_5^4}+3920000000q^6\dfrac{E_{2}^{8}E_{10}^{18}}{E_1^{20}E_5^5}\notag\\
&\quad+7680000000q^7\dfrac{E_{2}^{9}E_{10}^{21}}{E_1^{23}E_5^6}+6400000000q^8\dfrac{E_{2}^{10}E_{10}^{24}}{E_1^{26}E_5^7}.
\end{align*}
Using the above in \eqref{omega-main}, we readily arrive at \eqref{gen-om-40n12-main} to finish the proof.
\end{pf-thm-p-omega}

\begin{pf-cor-p-omega}
Taking congruences modulo 5 on both sides of \eqref{gen-om-40n12-main}, we have
\begin{align}\label{gen-om-40n12-cor1}\sum_{n=0}^\infty p_\omega(40n+12)q^n&\equiv\dfrac{E_{2}^{2}E_5}{E_1^{2}}~(\textup{mod}~5).\end{align}
Employing \eqref{E1} and \eqref{1byE1} in the above, extracting the terms involving $q^{5n+2}$, dividing both sides by $q^2$, and then replacing $q^5$ by $q$, we obtain
\begin{align*}\sum_{n=0}^\infty p_\omega(200n+92)q^n&\equiv3\dfrac{E_{2}^{2}E_5^2}{E_1^{3}}\equiv3E_1^{2}E_{2}^{2}E_5~(\textup{mod}~5),\end{align*}
where the last congruence is by \eqref{binomial}. Once again invoking \eqref{E1} in the above, extracting the terms involving $q^{5n+1}$, dividing both sides by $q$, and then replacing $q^5$ by $q$, we find that
\begin{align}\label{gen-om-40n12-cor2}\sum_{n=0}^\infty p_\omega(1000n+292)q^n&\equiv-\dfrac{E_1^{8}E_{10}}{E_2^3E_5}\equiv-\dfrac{E_2^{2}E_{5}}{E_1^2}~(\textup{mod}~5).\end{align}From \eqref{gen-om-40n12-cor1} and \eqref{gen-om-40n12-cor2}, we see that
\begin{align*}p_\omega(1000n+292)=p_\omega(40(25n+7)+12)&\equiv - p_\omega(40n+12)~(\textup{mod}~5).\end{align*}Iterating the above congruence as shown by Xia \cite{xia-mock}, one can easily arrive at \eqref{p-omega-inf} to finish the proof.

\end{pf-cor-p-omega}

\section{\textbf{Proofs of Theorem \ref{thm-spt-main} and the congruences \eqref{cong250n73} and \eqref{cong1250n573}}}\label{sec5}

\begin{pf-thm-spt-main}
Employing \eqref{A-5qB} in \eqref{gen-spt-om-2n1}, we find that
\begin{align}
\sum_{n=0}^\infty\textup{spt}_\omega\left(2n+1\right)q^n=\dfrac{E_2^3E_{10}}{E_5}+5q\dfrac{E_1^2E_{10}^5}{E_2E_5^3}+25q^2\dfrac{E_{10}^8}{E_1E_5^4}
\end{align}
Applying \eqref{E1} and \eqref{1byE1} in the above, extracting the terms involving $q^{5n+1}$, dividing both sides by $q$, replacing $q^5$ by $q$, and then proceeding as in the previous section, we obtain \eqref{gen-spt-10n3-main}.

With the aid of \eqref{A-5qB}, we can rewrite \eqref{gen-spt-10n3-main} as
\begin{align}\label{gen-spt-10n3-main1}\sum_{n=0}^\infty\textup{spt}_\omega\left(10n+3\right)q^n
&=5\Bigg(E_1^2E_5+6q\dfrac{E_2E_{10}^3}{E_1}+25q\dfrac{E_5^3E_{10}^2}{E_2^{2}}+250q^2\dfrac{E_5^2E_{10}^5}{E_1^{3}E_2}\notag\\
&\quad+
625q^3\dfrac{E_5E_{10}^8}{E_1^6}\Bigg).
\end{align}

Now, let $\left[q^{5n+r}\right]\left\{F(q)\right\}$, $r=0,1,\ldots, 4$ denote
 the terms after extracting the terms involving $q^{5n+r}$, dividing by $q^r$ and then replacing $q^5$ by $q$.

With the aid of Lemmas \ref{dissect}--\ref{AB}, omitting details, we find that
\begin{align*}
&\left[q^{5n+2}\right]\left\{E_1^2E_5+6q\dfrac{E_2E_{10}^3}{E_1}\right\}=5\left(E_1E_5^2+6q\dfrac{E_2E_5E_{10}^{3}}{E_1^{2}}\right),\\
&\left[q^{5n+2}\right]\left\{25q\dfrac{E_5^3E_{10}^2}{E_2^{2}}\right\}=125\left(2q\dfrac{E_1^3E_{10}^4}{E_2^4}+25q^3\dfrac{E_1^3E_{10}^{10}}{E_2^{10}}\right),\\
&\left[q^{5n+2}\right]\left\{250q^2\dfrac{E_5^2E_{10}^5}{E_1^{3}E_2}\right\}\\
&=250\Bigg(\dfrac{E_1^{16}E_{10}^{3}}{E_2^{15}E_5}+150q
\dfrac{E_1^{13}E_{10}^{6}}{E_2^{14}E_5^2}+5650q^2\dfrac{E_1^{10}E_{10}^{9}}{E_2^{13}E_5^3}+101825q^3\dfrac{E_1^{7}E_{10}^{12}}{E_2^{12}E_5^4}\\
&\quad+1068125q^4\dfrac{E_1^{4}E_{10}^{15}}{E_2^{11}E_5^5}+7042500q^5\dfrac{E_1E_{10}^{18}}{E_2^{10}E_5^6}
+29800000q^6\dfrac{E_{10}^{21}}{E_1^{2}E_2^{9}E_5^7}\end{align*}
\begin{align*}
&\quad+79000000q^7\dfrac{E_{10}^{24}}{E_1^{5}E_2^{8}E_5^8}+120000000q^8\dfrac{E_{10}^{27}}{E_1^{8}E_2^{7}E_5^9}
+80000000q^9\dfrac{E_{10}^{30}}{E_1^{11}E_2^{6}E_5^{10}}\Bigg),\\
&\left[q^{5n+2}\right]\left\{625q^3\dfrac{E_5E_{10}^8}{E_1^6}\right\}\\
&=3125q\Bigg(63\dfrac{E_2^8E_5^6}{E_1^{11}}+6500q\dfrac{E_2^8E_5^{12}}{E_1^{17}}+196875q^2\dfrac{E_2^8E_5^{18}}{E_1^{23}}+2343750q^3\dfrac{E_2^8E_5^{24}}{E_1^{29}}\notag\\
&\quad+9765625q^4\dfrac{E_2^8E_5^{30}}{E_1^{35}}\Bigg).
\end{align*}
Invoking the above in \eqref{gen-spt-10n3-main1}, we obtain \eqref{gen-spt-50n23-main}, as desired.
\end{pf-thm-spt-main}

\begin{pf-spt-1250n573}
Taking congruences modulo 625 on both sides of \eqref{gen-spt-50n23-main}, we have
\begin{align*}\sum_{n=0}^\infty\textup{spt}_\omega\left(50n+23\right)q^n
&\equiv25\left(E_1E_5^2+6q\dfrac{E_2E_5E_{10}^3}{E_1^2}\right)~(\textup{mod}~625).
\end{align*}
Employing \eqref{E1} and \eqref{1byE1} in the above, extracting the terms involving $q^{5n+1}$, dividing both sides by $q$, replacing $q^5$ by $q$ and then proceeding as in the earlier sections, we obtain
\begin{align*}\sum_{n=0}^\infty\textup{spt}_\omega\left(250n+73\right)q^n&\equiv 25\Bigg(-E_1^2E_5+6\dfrac{E_2^5E_5^2}{E_1^3E_{10}}+120q\dfrac{E_2^6E_5E_{10}^2}{E_1^6}\\
&\quad+480q^2\dfrac{E_2^7E_{10}^5}{E_1^9}\Bigg)~(\textup{mod}~625).
\end{align*}
Employing \eqref{A-5qB} and \eqref{binomial} in the above, we have
\begin{align*}&\sum_{n=0}^\infty\textup{spt}_\omega\left(250n+73\right)q^n\\
&\equiv 25\Bigg(-E_1^2E_5+6\left(E_1^2E_5+5q\dfrac{E_2E_{10}^3}{E_1}\right)+120q\left(\dfrac{E_2E_{10}^3}{E_1}+5q\dfrac{E_2^2E_{10}^6}{E_1^4E_5}\right)\\
&\quad+480q^2\left(\dfrac{E_2^2E_{10}^6}{E_1^4E_5}+5q\dfrac{E_2^3E_{10}^9}{E_1^7E_5^2}\right)\Bigg)\\
&\equiv125\left(E_1^2E_5+ 216q^2\dfrac{E_2^2E_{10}^6}{E_1^4E_5}\right)\\
&\equiv125\left(E_1^2E_5+q^2\dfrac{E_1E_2^2E_{10}^6}{E_{5}^2}\right)~(\textup{mod}~625),
\end{align*}
which readily implies \eqref{cong250n73}. Again employing \eqref{E1} in the above, extracting the terms involving $q^{5n+2}$, dividing both sides by $q^2$, replacing $q^5$ by $q$, and then proceeding as before, we find that
\begin{align}\label{spt-w-11}\sum_{n=0}^\infty\textup{spt}_\omega\left(1250n+573\right)q^n&\equiv125\left(-E_1E_5^2+\dfrac{E_2^{10}E_{5}^3}{E_1^4E_{10}^2}\right)~(\textup{mod}~625).
\end{align}
Since by \eqref{binomial},
\begin{align*}-E_1E_5^2+\dfrac{E_2^{10}E_{5}^3}{E_1^4E_{10}^2}&\equiv0~(\textup{mod}~5),
\end{align*}
from \eqref{spt-w-11}, we readily arrive at \eqref{cong1250n573} to finish the proof.
\end{pf-spt-1250n573}

\section{\textbf{Proofs of Theorem \ref{thm-sptbar} and Corollary \ref{cor-sptbar}}}\label{sec6}

\begin{pf-thm-sptbar}Employing \eqref{A-4qB} -- \eqref{A-5qB} in \eqref{gen-overspt-om-2n1}, we see that
\begin{align*}\sum_{n=0}^\infty \overline{\textup{spt}}_\omega\left(2n+1\right)q^n&=\dfrac{E_2^4E_{10}}{E_1E_5}+5q\dfrac{E_2^5E_{10}^4}{E_1^4E_5^2}\notag\\
&=E_2E_5^2+q\dfrac{E_1E_{10}^{5}}{E_5^{3}}+5q\left(\dfrac{E_1E_{10}^{5}}{E_5^{3}}+5q\dfrac{E_2E_{10}^{8}}{E_1^2E_5^{4}}\right)\notag\\
&=E_2E_5^2+6q\dfrac{E_1E_{10}^{5}}{E_5^{3}}+25q^2\left(\dfrac{E_1^{2}E_{10}^{10}}{E_2E_5^{8}}
+4q\dfrac{E_{10}^{13}}{E_1E_{5}^9}\right)\\
&=E_2E_5^2+6q\dfrac{E_1E_{10}^{5}}{E_5^{3}}+25q^2\dfrac{E_1^{2}E_{10}^{10}}{E_2E_5^{8}}
+100q^3\dfrac{E_{10}^{13}}{E_1E_{5}^9}.
\end{align*}
Applying \eqref{E1} and \eqref{1byE1} in the above, extracting the terms involving $q^{5n+2}$, dividing both sides by $q^2$, replacing $q^5$ by $q$, we find that
\begin{align*}\sum_{n=0}^\infty \overline{\textup{spt}}_\omega\left(10n+5\right)q^n&=-E_1^2E_{10}-6\dfrac{E_2^5E_5}{E_1^3}+25\dfrac{E_2^{10}E_5^2}{E_1^{8}E_{10}}
+500q\dfrac{E_2^{13}E_5^5}{E_1^{15}}.
\end{align*}
We now apply \eqref{A-4qB} and \eqref{A-5qB} in the above to arrive at \eqref{gen-over-spt-10n5-main} as shown below:
\begin{align*}&\sum_{n=0}^\infty \overline{\textup{spt}}_\omega\left(10n+5\right)q^n\\
&=-E_1^2E_{10}-6\left(E_1^2E_{10}+5q\dfrac{E_2E_{10}^4}{E_1E_5}\right)
+25\left(\dfrac{E_2^{5}E_5}{E_1^{3}}+5q\dfrac{E_2^{6}E_{10}^3}{E_1^{6}}\right)\\
&\quad+500q\left(\dfrac{E_2^{8}E_5^4E_{10}}{E_1^{10}}+5q\dfrac{E_2^{9}E_5^3E_{10}^4}{E_1^{13}}\right)\\
&=-7E_1^2E_{10}-30q\dfrac{E_2E_{10}^4}{E_1E_5}
+25\left(E_1^2E_{10}+5q\dfrac{E_2E_{10}^4}{E_1E_5}\right)+125q\left(\dfrac{E_2E_{10}^4}{E_1E_5}
+5q\dfrac{E_2^2E_{10}^7}{E_1^4E_5^2}\right)\\
&\quad+500q\left(\dfrac{E_2^{3}E_5^3E_{10}^2}{E_1^{5}}
+5q\dfrac{E_2^{4}E_5^2E_{10}^5}{E_1^{8}}\right)+2500q^2\left(\dfrac{E_2^{4}E_5^2E_{10}^5}{E_1^{8}}
+5q\dfrac{E_2^{5}E_5E_{10}^8}{E_1^{11}}\right)\\
&=18E_1^2E_{10}+220q\dfrac{E_2E_{10}^4}{E_1E_5}+625q^2\dfrac{E_2^2E_{10}^7}{E_1^4E_5^2}+500q\dfrac{E_2^3E_5^3E_{10}^2}{E_1^5}+5000q^2\dfrac{E_2^{4}E_5^2E_{10}^5}{E_1^{8}}\\
&\quad+12500q^3\dfrac{E_2^{5}E_5E_{10}^8}{E_1^{11}}\\
&=18E_1^2E_{10}+220q\dfrac{E_2E_{10}^4}{E_1E_5}+625q^2\dfrac{E_2^2E_{10}^7}{E_1^4E_5^2}
+500q\Bigg(\dfrac{E_2E_{10}^4}{E_1E_5}+4q\dfrac{E_2^2E_{10}^7}{E_1^4E_5^2}\Bigg)\\
&\quad+5000q^2\Bigg(\dfrac{E_2^2E_{10}^7}{E_1^4E_5^2}+4q\dfrac{E_2^{3}E_{10}^{10}}{E_1^{7}E_5^3}\Bigg)
+12500q^3\Bigg(\dfrac{E_2^{3}E_{10}^{10}}{E_1^{7}E_5^3}+4q\dfrac{E_2^{4}E_{10}^{13}}{E_1^{10}E_5^4}\Bigg)\end{align*}\begin{align*}
&=18E_1^2E_{10}+720q\dfrac{E_2E_{10}^4}{E_1E_5}+7625q^2
\dfrac{E_2^2E_{10}^{7}}{E_1^{4}E_5^2}+32500q^3\dfrac{E_{2}^{3}E_{10}^{10}}{E_1^{7}E_5^3}+50000q^4\dfrac{E_{2}^{4}E_{10}^{13}}{E_1^{10}E_5^4}.
\end{align*}

The identity \eqref{gen-over-spt-50n25-main} can be proved in a similar way, therefore we omit the proof.

\end{pf-thm-sptbar}

\begin{pf-cor-sptbar} From \eqref{gen-over-spt-10n5-main} and \eqref{gen-over-spt-50n25-main}, we have
\begin{align}
\label{gen-over-spt-10n5-cor1}\sum_{n=0}^\infty \overline{\textup{spt}}_\omega\left(10n+5\right)q^n&\equiv 3E_1^2E_{10}~(\textup{mod}~5) \\\intertext{and}
\label{gen-over-spt-10n5-cor2}\sum_{n=0}^\infty \overline{\textup{spt}}_\omega\left(50n+25\right)q^n&\equiv 2E_2E_5^2~(\textup{mod}~5),
\end{align}
respectively. Employing \eqref{E1} in \eqref{gen-over-spt-10n5-cor2}, then extracting the terms involving $q^{5n+2}$, we arrive at
 \begin{align}
 \label{gen-over-spt-10n5-cor5}\sum_{n=0}^\infty \overline{\textup{spt}}_\omega\left(250n+125\right)q^n&\equiv 3E_1^2E_{10}~(\textup{mod}~5).
\end{align}
It follows from \eqref{gen-over-spt-10n5-cor1} and \eqref{gen-over-spt-10n5-cor5} that
\begin{align*}
  \overline{\textup{spt}}_\omega\left(250n+125\right)&\equiv \overline{\textup{spt}}_\omega\left(10n+5\right)~(\textup{mod}~5).
\end{align*}
Iterating the above, we easily deduce the case $\ell=1$ of \eqref{spt-gencong}.
%%%%%%%%%%%%%%%%%%%%%%%%%%%%%%%%%%%%%%%%%%%%%%%%%%%%%%%%%%%%%%%%%%%%%%%%%%%%%%%%%%%%%%%%%%%%%%%%%%%%%%%%%%%%%%%%%%%%%%%%%%%%%%%%%

Next, from \eqref{gen-over-spt-50n25-main} and \eqref{E1}, we have
\begin{align}
\label{25a}\sum_{n=0}^\infty \overline{\textup{spt}}_\omega\left(50n+25\right)q^n&\equiv2E_2E_5^2\notag\\
&\equiv2E_5^2E_{50}\left(T(q^{10})-q^2-\dfrac{q^{4}}{T(q^{10})}\right)~(\textup{mod}~5^2).
\end{align}
Extracting the terms involving $q^{5n+1}$ and $q^{5n+3}$, we arrive at \eqref{spt-cor2-mod5}. On the other hand,  extraction of the terms involving $q^{5n+2}$ gives
 \begin{align*}
\sum_{n=0}^\infty \overline{\textup{spt}}_\omega\left(250n+125\right)q^n&\equiv -2E_1^2E_{10}\notag\\
&\equiv -2E_{10}E_{25}^2\left(T(q^{5})-q-\dfrac{q^{2}}{T(q^{5})}\right)^2~(\textup{mod}~5^2),
\end{align*}
which implies
\begin{align}
\label{25e}\sum_{n=0}^\infty \overline{\textup{spt}}_\omega\left(25(50n+25)\right)q^n&\equiv 2E_2E_{5}^2~(\textup{mod}~5^2).
\end{align}
From \eqref{25a} and \eqref{25e}, we have
\begin{align*}
  \overline{\textup{spt}}_\omega\left(25(50n+25)\right)&\equiv \overline{\textup{spt}}_\omega\left(50n+25\right)~(\textup{mod}~5^2),
\end{align*}
which upon iteration proves the case $\ell=2$ of \eqref{spt-gencong}.
%%%%%%%%%%%%%%%%%%%%%%%%%%%%%%%%%%%%%%%%%%%%%%%%%%%%%%%%%%%%%%%%%%%%%%%%%%%%%%%%%%%%%%%%%%%%%%%%%%%%%%%%%%%%%%%%%%%%%%%%%%%%%%%%%%%%%%%%%%%%%%%

Now, from \eqref{gen-over-spt-50n25-main}, we also have
\begin{align*}
\sum_{n=0}^\infty \overline{\textup{spt}}_\omega\left(50n+25\right)q^n\equiv77E_2E_5^2+100q\dfrac{E_2^2E_5E_{10}^3}{E_1^3}~(\textup{mod}~5^3).
\end{align*}
Employing \eqref{E1} and \eqref{1byE1} in the above, extracting the terms involving $q^{5n+2}$, and then using Lemma \ref{lemmax-y}, we obtain
 \begin{align}
\label{125b}\sum_{n=0}^\infty \overline{\textup{spt}}_\omega\left(250n+125\right)q^n&\equiv 98E_1^2E_{10}\notag\\
&\equiv 98E_{10}E_{25}^2\left(T(q^{5})-q-\dfrac{q^{2}}{T(q^{5})}\right)^2~(\textup{mod}~5^3).
\end{align}
Therefore,
\begin{align*}
\sum_{n=0}^\infty \overline{\textup{spt}}_\omega\left(1250n+625\right)q^n&\equiv -98E_2E_{5}^2\notag\\
&\equiv -98E_{5}^2E_{50}\left(T(q^{10})-q^2-\dfrac{q^{4}}{T(q^{10})}\right)~(\textup{mod}~5^3).
\end{align*}
Extracting the terms involving $q^{5n+2}$, we arrive at 
 \begin{align*}
\sum_{n=0}^\infty \overline{\textup{spt}}_\omega\left(6250n+3125\right)q^n&\equiv 98E_1^2E_{10}~(\textup{mod}~5^3).
\end{align*}
It follows from the above and \eqref{125b} that
\begin{align*}
  \overline{\textup{spt}}_\omega\left(25(250n+125)\right)&\equiv \overline{\textup{spt}}_\omega\left(250n+125\right)~(\textup{mod}~5^3).
\end{align*}
Iterating the above, we easily deduce the case $\ell=3$ of \eqref{spt-gencong}.

The remaining cases of \eqref{spt-gencong} can be proved in a similar fashion. Therefore, we omit the details and record only the successive generating functions.

For congruences modulo $5^4$, we have
\begin{align}
\sum_{n=0}^\infty \overline{\textup{spt}}_\omega\left(50n+25\right)q^n&\equiv202E_2E_5^2+475q\dfrac{E_2^2E_5E_{10}^3}{E_1^3},\notag\\
\sum_{n=0}^\infty \overline{\textup{spt}}_\omega\left(250n+125\right)q^n&\equiv 598E_1^2E_{10}+250q\dfrac{E_2E_{10}^4}{E_1E_5}+125q^2\dfrac{E_2^2E_{10}^7}{E_1^4E_5^2},\notag\\
\label{5power4}\sum_{n=0}^\infty \overline{\textup{spt}}_\omega\left(1250n+625\right)q^n&\equiv -223E_2E_{5}^2,\\
\sum_{n=0}^\infty \overline{\textup{spt}}_\omega\left(6250n+3125\right)q^n&\equiv 223E_1^2E_{10},\notag\\\intertext{and}
\sum_{n=0}^\infty \overline{\textup{spt}}_\omega\left(31250n+15625\right)q^n&\equiv -223E_2E_{5}^2.\notag
\end{align}

For congruences modulo $5^5$, we have

\begin{align*}
\sum_{n=0}^\infty \overline{\textup{spt}}_\omega\left(50n+25\right)q^n&\equiv2077E_2E_5^2+2975q\dfrac{E_2^2E_5E_{10}^3}{E_1^3},\\
\sum_{n=0}^\infty \overline{\textup{spt}}_\omega\left(250n+125\right)q^n&\equiv 598E_1^2E_{10}+1500q\dfrac{E_2E_{10}^4}{E_1E_5}+1375q^2\dfrac{E_2^2E_{10}^7}{E_1^4E_5^2}\\
&\quad+1250q^3\dfrac{E_2^3E_{10}^{10}}{E_1^2E_5^4},\\
\sum_{n=0}^\infty \overline{\textup{spt}}_\omega\left(1250n+625)\right)q^n&\equiv2277E_2E_{5}^2 +625q\dfrac{E_2^2E_5E_{10}^3}{E_1^3},\\
\sum_{n=0}^\infty \overline{\textup{spt}}_\omega\left(6250n+3125\right)q^n&\equiv -402E_1^2E_{10},\\
\sum_{n=0}^\infty \overline{\textup{spt}}_\omega\left(31250n+15625\right)q^n&\equiv 402E_2E_{5}^2,\\\intertext{and}
\sum_{n=0}^\infty \overline{\textup{spt}}_\omega\left(156250n+78125\right)q^n&\equiv -402E_1^2E_{10}.
\end{align*}

Finally, for congruences modulo $5^6$, we have

\begin{align}
\sum_{n=0}^\infty \overline{\textup{spt}}_\omega\left(50n+25\right)q^n&\equiv8327E_2E_5^2+15475q\dfrac{E_2^2E_5E_{10}^3}{E_1^3}+6250q^2\dfrac{E_2^3E_{10}^6}{E_1E_5},\notag\\
\sum_{n=0}^\infty \overline{\textup{spt}}_\omega\left(250n+125\right)q^n&\equiv 13098E_1^2E_{10}+4625q\dfrac{E_2E_{10}^4}{E_1E_5}+4500q^2\dfrac{E_2^2E_{10}^7}{E_1^4E_5^2}\notag\\
&\quad+1250q^3\dfrac{E_2^3E_{10}^{10}}{E_1^7E_5^3}+3125q^4\dfrac{E_{10}^{14}}{E_2E_5^6},\notag\\
\sum_{n=0}^\infty \overline{\textup{spt}}_\omega\left(1250n+625)\right)q^n&\equiv-3973E_2E_{5}^2 +6875q\dfrac{E_2^2E_5E_{10}^3}{E_1^3}+3125q^2\dfrac{E_2^3E_{10}^6}{E_1E_5},\notag\\
\sum_{n=0}^\infty \overline{\textup{spt}}_\omega\left(6250n+3125\right)q^n&\equiv 12098E_1^2E_{10}+12500q\dfrac{E_2E_{10}^4}{E_1E_5}+12500q^2\dfrac{E_2^2E_{10}^7}{E_1^4E_5^2},\notag\\
\label{5power6}\sum_{n=0}^\infty \overline{\textup{spt}}_\omega\left(31250n+15625\right)q^n&\equiv 12902E_2E_{5}^2,\\
\sum_{n=0}^\infty \overline{\textup{spt}}_\omega\left(156250n+78125\right)q^n&\equiv -12902E_1^2E_{10},\notag\\\intertext{and}
\sum_{n=0}^\infty \overline{\textup{spt}}_\omega\left(781250n+390625\right)q^n&\equiv 12902E_2E_{5}^2.\notag
\end{align}
Employing \eqref{E1}, with $q$ replaced by $q^2$, in \eqref{5power4} and \eqref{5power6}, and then extracting the terms involving $q^{5n+1}$ and $q^{5n+3}$ from the resulting congruential identities, we arrive at \eqref{spt-cor2-mod125} and \eqref{spt-cor2-mod3125} to finish the proof.
\end{pf-cor-sptbar}

\section*{Acknowledgment} The authors would like to thank the referee for reading the manuscript with meticulous care, uncovering several errors and offering his/her helpful suggestions. The first author's research was partially supported by Grant no. MTR/2018/000157 of Science \& Engineering Research Board (SERB), DST, Government of India.

\end{document}